\documentclass[a4,12pt]{amsart}
\usepackage{amssymb,amstext,amsmath,amscd,amsthm,amsfonts,enumerate,latexsym, comment}
\usepackage{color}
\usepackage[all]{xy}

\theoremstyle{plain}
\newtheorem{theorem}{Theorem}[section]
\newtheorem{proposition}[theorem]{Proposition}
\newtheorem{corollary}[theorem]{Corollary}
\newtheorem{lemma}[theorem]{Lemma}

\theoremstyle{definition}
\newtheorem{setting}[theorem]{Setting}
\newtheorem{rem}[theorem]{Remark}
\newtheorem{fact}[theorem]{Fact}
\newtheorem{ex}[theorem]{Example}
\newtheorem*{acknowledgments}{Acknowledgments}

\def\Min{\operatorname{Min}}

\def\rank{\mathrm{rank}}

\newcommand{\rme}{\mathrm{e}}

\newcommand{\calC}{\mathcal{C}}

\newcommand{\calG}{\mathcal{G}}

\newcommand{\calL}{\mathcal{L}}
\newcommand{\calM}{\mathcal{M}}

\newcommand{\calR}{\mathcal{R}}
\newcommand{\calS}{\mathcal{S}}

\newcommand{\fka}{\mathfrak{a}}

\newcommand{\fkm}{\mathfrak{m}}

\newcommand{\fkp}{\mathfrak{p}}
\newcommand{\fkq}{\mathfrak{q}}

\newcommand{\mapright}[1]{%
\smash{\mathop{%
\hbox to 1cm{\rightarrowfill}}\limits^{#1}}}

\newcommand{\mapleft}[1]{%
\smash{\mathop{%
\hbox to 1cm{\leftarrowfill}}\limits_{#1}}}

\def\depth{\operatorname{depth}}

\def\Ass{\operatorname{Ass}}
\def\height{\mathrm{ht}}

\def\Spec{\operatorname{Spec}}

\def\ol{\overline}

\title{graded filtrations and ideals of reduction number two}
\author{Shinya Kumashiro}
\address{National Institute of Technology (KOSEN), Oyama College, 771 Nakakuki, Oyama, Tochigi, 323-0806, Japan}
\email{skumashiro@oyama-ct.ac.jp}


\begin{document}

\begin{abstract} 
In this paper, we give a way to construct graded filtrations of graded modules. We then apply it to the Sally module, which describes a correction term of the Hilbert function. As a result, we obtain the inequality of the Hilbert coefficients for ideals of reduction number $2$ or $3$.
\end{abstract}

\maketitle

\section{Introduction}\label{section1}

In this paper, we give a new construction of graded filtration that can control the shifts for finitely generated graded modules. By using this new filtration, we later obtain the inequality of the first three Hilbert coefficients. Let us present our notation to state our results precisely. Let $R=\bigoplus_{n\ge 0}R_n$ be a standard graded Noetherian ring of dimension $\ge 2$ such that $R_0$ is an Artinian local ring with the unique maximal ideal $\fkm_0$. Recall that for a finitely generated graded $R$-module $M$ of dimension $s$, there exist the integers $\rme_0(M), \rme_1(M), \dots, \rme_{s-1}(M)$ such that 
\[
\ell_{R_0}(M_n)=\rme_0(M) \binom{n+s-1}{s-1} - \rme_1(M) \binom{n+s-2}{s-2} + \cdots +(-1)^{s-1}\rme_{s-1}(M)
\]
for $n \gg 0$ (\cite[Theorem 4.1.3]{BH}), where $\ell(*)$ denotes the length. $\rme_0(M), \rme_1(M), \dots, \rme_{s-1}(M)$ are called the {\it Hilbert coefficients} of $M$. Then, the first main theorem of this paper can be stated as follows:

\begin{theorem} \label{thm1}
Let $R=\bigoplus_{n\ge 0}R_n$ be a standard graded Noetherian ring of dimension $\ge 2$ such that $R_0$ is an Artinian local ring with the unique maximal ideal $\fkm_0$. Suppose that $\fkm_0 R$ is a prime ideal, and set $\fkp=\fkm_0 R$. 
Let $M$ be a finitely generated graded $R$-module generated in single degree $t$.  Suppose that $\dim_R M=\dim R$ and  $\Ass_R M=\Min_R M$. Then, the following hold true.
\begin{enumerate}[{\rm (a)}] 
\item $\rme_0(M)=\ell_{R_\fkp}(M_\fkp){\cdot}\rme_0(R/\fkp)$\quad and \quad $\rme_1(M)\le t\rme_0(M) + \ell_{R_\fkp}(M_\fkp){\cdot}\rme_1(R/\fkp)$.
\item Consider the following conditions:
\begin{enumerate}[{\rm (i)}] 
\item $\rme_1(M)= t\rme_0(M) + \ell_{R_\fkp}(M_\fkp){\cdot}\rme_1(R/\fkp)$.
\item There exists a family $\{M^i\}_{0\le i\le i_0}$ of graded residue modules of $M$ such that
\begin{align}\label{**} 
\begin{split} 
0 \to (R/\fkp)(-t) \to &M=M^0 \to M^1 \to 0,\\
0 \to (R/\fkp)(-t) \to &M^1 \to M^2 \to 0,\\
&\vdots\\
0 \to (R/\fkp)(-t) \to &M^{i_0-2} \to M^{i_0-1} \to 0, \text{ and}\\
0 \to (R/\fkp)(-t) \to &M^{i_0-1} \to M^{i_0}=0 \to 0
\end{split}
\end{align}
are exact.
\end{enumerate} 
Then, {\rm (ii) $\Rightarrow$ (i)} holds. {\rm (i) $\Rightarrow$ (ii)} also holds if $R/\fkp$ satisfies Serre's condition $(S_2)$. 
\end{enumerate} 
\end{theorem}

We can further provide an application of Theorem \ref{thm1} to the Hilbert function of the associated graded rings. Let $(A, \fkm)$ be a Cohen--Macaulay local ring of dimension $d\ge 2$, and let  $I$ be an $\fkm$-primary ideal of $A$. Let 
\begin{align*} 
\calR(I)=A[It] \subseteq A[t] \quad \text{and} \quad  \calG(I)=\calR(I)/I\calR(I)\cong \bigoplus_{n\ge 0}I^n/I^{n+1}
\end{align*}
be the {\it Rees algebra} of $I$ and the {\it associated graded ring} of $I$, respectively, where $A[t]$ is the polynomial ring over $A$. 
Then, there exist the integers $\rme_0(I), \rme_1(I), \dots, \rme_{d}(I)$ such that 
\[
\ell_A(A/I^{n+1})=\rme_0(I) \binom{n+d}{d} - \rme_1(I) \binom{n+d-1}{d-1} + \cdots +(-1)^{d}\rme_{d}(I)
\]
for $n \gg 0$ because $\ell_A(A/I^{n+1}) - \ell_A(A/I^{n}) = \ell_A(\calG(I)_n)$ for all $n \ge 0$. $\rme_0(I), \rme_1(I), \dots, \rme_{d}(I)$ are called the {\it Hilbert coefficients} of $I$. The problem to reveal the relation between the Hilbert coefficients of $I$ and the structures of $\calR(I)$ and $\calG(I)$ is studied very much by many researchers (see papers \cite{CPR, GNO, Hu, H, I, K2, K, Ma1, No, O, OR, S2, V, VP} or a book \cite{RV}). One of the most famous results in this direction is the following obtained by Northcott, Huneke, and Ooishi. 

\begin{fact}\label{Northcott} {\rm (\cite{No, H, O})}
Let $(A, \fkm)$ be a Cohen--Macaulay local ring of dimension $d\ge 2$, and let  $I$ be an $\fkm$-primary ideal of $A$. Suppose that $A/\fkm$ is an infinite field. Then 
\[
\ell_A(A/I)\ge \rme_0(I) - \rme_1(I)
\] 
holds. The equality holds if and only if $I^2=QI$ for some (any) minimal reduction $Q\subseteq I$. When the equality holds, the Rees algebra of $I$ and the associated graded ring of $I$ are Cohen--Macaulay rings.
\end{fact}

The second aim of this paper is to study the next step, i.e., the case of reduction number $2$. Note that the {\it reduction number} (with respect to $Q$), the minimal integer $n\ge 0$ satisfying $I^{n+1}=QI^n$, depends on the choice of a minimal reduction $Q$ in general (see for example \cite{Hu, Ma1, Ma2}); hence, the behavior of the Hilbert function is complicated even if the case of reduction number $2$. For this progress, see for example \cite[Introduction]{K}. The following theorem clarifies what happen in the case of reduction number $2$. 

\begin{theorem}{\rm (Theorem \ref{thm3.3})}\label{thm2}
In addition to the assumption of {\rm Fact \ref{Northcott}}, suppose that $I^3=QI^2$ for some minimal reduction $Q\subseteq I$. Then, 
\[
\ell_A(A/I)\ge \rme_0(I) - \rme_1(I) + \rme_2(I)
\] 
holds. The equality holds if and only if $\depth \calG(I) \ge d-1$.
\end{theorem}

How to apply Theorem \ref{thm1} to Theorem \ref{thm2} is a bit complicated. We will apply Theorem \ref{thm1} to the Sally module, which describes a correction term of the Hilbert function. The Sally module can be regarded as a finitely generated graded $\calR(Q)/\fkm^\ell \calR(Q)$-module for $\ell \gg 0$, see the observation before Proposition \ref{prop3.2}. 


In contrast Theorem \ref{thm2}, it is known that the reverse inequality 
\[
\ell_A(A/I)\le \rme_0(I) - \rme_1(I) + \rme_2(I)
\] 
holds if $I$ is either an integrally closed ideal or a Ratliff--Rush closed ideal of height two (\cite{S2, I}).  Thus, Theorem \ref{thm2} also provides a deeper understanding of a result of Corso-Polini-Rossi (\cite[Theorem 3.6]{CPR}).

The remainder of this paper is organized as follows. In Section \ref{section2}, we construct a graded filtration and establish Theorem \ref{thm1}. In Section \ref{section3}, we demonstrate the application of Theorem \ref{thm1} to the Sally module and prove Theorem \ref{thm2}. We further examine the case of reduction number $3$ under the assumption that $I$ is integrally closed (Theorem \ref{final}).


\section{Graded filtration of graded modules}\label{section2}
The aim of this section is to construct a new graded filtration. 
We maintain the following setting throughout this section.

\begin{setting} 
Let $R=\bigoplus_{n\ge 0}R_n=R_0[R_1]$ be a standard graded Noetherian ring such that $R_0$ is an Artinian local ring with the maximal ideal $\fkm_0$. Let $M=\bigoplus_{n\in \mathbb{Z}}M_n$ be a finitely generated graded $R$-module. Suppose that $\fkm_0 R$ is a prime ideal, and set $\fkp=\fkm_0 R$.
\end{setting}

\begin{rem} 
\begin{enumerate}[{\rm (a)}] 
\item $\fkm_0^\ell M=0$ for $\ell \gg 0$.
\item $\Min R=\{\fkp\}$.
\item The following are equivalent:
\begin{enumerate}[{\rm (i)}] 
\item $\dim_R M=\dim R$ and $\Ass_R M=\Min_R M$.
\item $\Ass_R M=\{\fkp\}$.
\end{enumerate} 
\end{enumerate} 
\end{rem}

\begin{proof} 
(a):  Because $(R_0, \fkm_0)$ is an Artinian local ring, $\ell_{R_0} (M_n)< \infty$. It follows that $\fkm_0^\ell M=0$ for  $\ell \gg 0$ because $M$ is finitely generated.

(b): One hand, we have that $\displaystyle \sqrt{(0)}=\bigcap_{\fkq\in \Spec R} \fkq\subseteq \fkp$, where $\sqrt{(0)}$ denotes the nilradical of $R$.
On the other hand, $\fkm_0^\ell R=0$ for $\ell \gg 0$ by (a); hence, $\fkp\subseteq \sqrt{(0)}$. It follows that $\displaystyle \fkp=\sqrt{(0)}=\bigcap_{\fkq\in \Spec R} \fkq$, i.e., $\Min R=\{\fkp\}$.

(c)(i)$\Rightarrow$(ii): Because $\dim_R M=\dim R$, $\fkp\in \Min_R M$. It follows that $\Ass_R M=\Min_R M=\{\fkp\}$ by (b). 

(c)(ii)$\Rightarrow$(i): This is clear.
\end{proof}

\begin{lemma} {\rm (cf. \cite[Chapter ${\rm I\hspace{-.15em}V}$, Section 1, 1. Proposition 4.]{B})}\label{2.3}
Let $\Phi$ be a subset of $\Ass_R M$. Then, there exists a graded exact sequence 
\[
0 \to L \to M \to N \to 0,
\]
where $L$ and $N$ are graded $R$-modules, such that $\Ass_R L=\Ass_R M\setminus \Phi$ and $\Ass_R N=\Phi$.
\end{lemma}

\begin{proof} 
This is a graded version of \cite[Chapter ${\rm I\hspace{-.15em}V}$, Section 1, 1. Proposition 4.]{B}, and the proof is essentially the same. However, we include the proof for the sake of completeness. 
We may assume that $\Phi\subsetneq \Ass_R M$. Then
\[
S=\left\{ X \ \middle| \  \begin{matrix}
\text{$X$ is a finitely generated graded $R$-submodule of M such that} \\
\text{$\Ass_R X \subseteq \Ass_R M \setminus \Phi$}
\end{matrix}
\right\}
\]
is a non-empty set. Hence, we can choose a maximal element $L$ with respect to the inclusion. 
We demonstrate that $\Ass_R M/L \subseteq \Phi$. Indeed, if this holds, then the assertion holds because $\Ass_R M \subseteq \Ass_R L \cup \Ass_R M/L$. Assume the contrary. Then, there exist $\fkq=(0):_R x\in \Ass_RM/L \setminus \Phi$, where $x$ is a homogeneous element of degree $t$, and an inclusion $0 \to (R/\fkq) (-t) \to M/L$ by \cite[Lemma 1.5.6(b)(ii)]{BH}. It follows that 
\[
\Ass_R (L+Rx)\subseteq \Ass_R M/L \setminus \Phi
\] 
by $0 \to L \to L+Rx \to (R/\fkq)(-t) \to 0$. This contradicts the maximality of $L$. Hence, $\Ass_R M/L \subseteq \Phi$.
\end{proof} 

\begin{lemma}\label{2.4}
Suppose that $\Ass_R M=\{\fkp\}$. If $M_t\ne 0$ for $t\in \mathbb{Z}$, then there exists a graded inclusion $0 \to (R/\fkp)(-t) \to M$.
\end{lemma}

\begin{proof} 
Because $R_0$ is an Artin local ring, $\fkm_0 \in \Ass_{R_0} M_t$. Thus, $\fkp=\fkm_0 R$ annihilates a non-zero element $x$ of $M_t$ as an element of the $R$-module $M$. Meanwhile, the annihilator of $x$ is contained in an associated prime ideal of $M$. Since $\fkp$ is the only prime ideal of $M$, the annihilator of $x$ is contained in $\fkp$. It follows that $\fkp$ is the annihilator of $x$. Hence, the $R$-homomorphism $R(-t) \to M; 1\mapsto x$ induces the required inclusion $0 \to (R/\fkp)(-t) \to M$.
\end{proof} 

The combination of Lemmas \ref{2.3} and \ref{2.4} provides a family of graded commutative diagrams as follows:

\begin{proposition}\label{filtration}
Suppose that $\Ass_R M=\{\fkp\}$ and $M=RM_t$. Then, there exist families of graded modules $\{M^i\}_{0\le i\le i_0}$ and $\{X^i\}_{0\le i\le i_0-1}$ such that $M^0=M$, $M^{i_0}=0$, and the following three assertions hold true for all $0\le i \le i_0-1$.
\begin{enumerate}[{\rm (a)}] 
\item $M^i$ satisfies the same assumption as $M$, i.e., $\Ass_R M^i=\{\fkp\}$ and $M^i=RM^i_t$, 
\item $\fkp \not\in \Ass_R X^i$, and 
\item there exists a graded commutative diagram 
\begin{align} \label{*}
\begin{split} 
\xymatrix{
&&&0 \ar[d]&\\
&&& X^i \ar[d] & \\
0 \ar[r] & (R/\fkp)(-t) \ar[r]  & M^i \ar[r] & N^i \ar[r] \ar[d] & 0 \\
&&&M^{i+1} \ar[d]&\\
&&&0&.\\
}
\end{split}
\end{align}
\end{enumerate} 
Furthermore, we then have $i_0=\ell_{R_{\fkp}}(M_\fkp)$.
\end{proposition}

\begin{proof} 
By Lemma \ref{2.4}, we have a graded exact sequence 
\[
0 \to (R/\fkp)(-t) \to M \xrightarrow{\pi_1} N^0 \to 0.
\] 
If $\fkp\not\in \Ass_R N^0$, then $M^1=0$ and $N^0=X^0$; hence, the assertion holds. Suppose that $\fkp\in \Ass_R N^0$. Then, by Lemma \ref{2.3}, we also have a graded exact sequence 
\[
0 \to X^0 \to N^0 \xrightarrow{\pi_2} M^1 \to 0
\] 
such that $\Ass_R X^0=\Ass_R N^0\setminus \{\fkp\}$ and $\Ass_R M^1= \{\fkp\}$. Because the composition $M \xrightarrow{\pi_1} N^0\xrightarrow{\pi_2} M^1$ is surjective, $M^1$ is also generated in single degree $t$. Hence, $M^1$ and $X^0$ satisfy the conditions (a) and (b), respectively, thus the two graded exact sequences above provide the graded commutative diagram \eqref{*} for $i=0$. Since $M^1$ satisfies the same assumption as $M$, this process can be continued recursively. Note that the process will end within a finite number of steps. Indeed, by looking at \eqref{*}, we obtain that $\ell_{R_0}(M^{i+1}_t)<\ell_{R_0}(M^i_t)$. It follows that $\ell_{R_0}(M^{i_0}_t)=0$ for $i_0\gg 0$; thus, $M^{i_0}=0$ since $M^{i_0}$ is generated in degree $t$.
Furthermore, by localizing (\ref{*}) at $\fkp$, we obtain that 
\[
\ell_{R_{\fkp}}(M^i_\fkp)=\ell_{R_{\fkp}}((R/\fkp)_\fkp) + \ell_{R_{\fkp}}(X^i_\fkp) + \ell_{R_{\fkp}}(M^{i+1}_\fkp) = 1 + \ell_{R_{\fkp}}(M^{i+1}_\fkp)
\] 
for all $0\le i \le i_0-1$. It follows that $\ell_{R_{\fkp}}(M_\fkp)=i_0$.
\end{proof}

Now we can prove Theorem \ref{thm1}, which characterizes when $X^i=0$ in  (\ref{*}) for all $0 \le i \le i_0-1$ in terms of the Hilbert coefficients.

\begin{proof}[Proof of Theorem \ref{thm1}]
In this proof, let $d=\dim R$. 

(a): By Proposition \ref{filtration}, we have that 
\[
\ell_{R_0}(M_n)=j_0{\cdot}\ell_{R_0}((R/\fkp)_{n-t}) + \sum_{i=0}^{j_0-1} \ell_{R_0}(X_n^i)
\]
for all $n\in \mathbb{Z}$, where $j_0=\ell_{R_\fkp}(M_\fkp)$. Because $\fkp\not\in \Ass_RX^i$, $\ell_{R_0}(X_n^i)$ is a polynomial of degree $<d-1$ for $n\gg 0$. Let $\rme(X^i)$ denote the coefficient of $\binom{n+d-2}{d-2}$ in the polynomial ($\rme(X^i)$ can be $0$). Then, 
\begin{align*} 
\ell_{R_0}(M_n)=&j_0\left(\rme_0(R/\fkp)\binom{n-t+d-1}{d-1} - \rme_1(R/\fkp)\binom{n-t+d-2}{d-2} + \cdots + (-1)^{d-1}\rme_{d-1}(R/\fkp)\right) \\
& + \sum_{i=0}^{j_0-1} \left(\rme(X^i)\binom{n+d-2}{d-2} + (\text{polynomial of degree $<d-2$})\right)\\
=&j_0{\cdot}\rme_0(R/\fkp) \binom{n+d-1}{d-1} + \left( j_0{\cdot}\rme_0(R/\fkp)(-t) - j_0{\cdot}\rme_1(R/\fkp) + \sum_{i=0}^{j_0-1} \rme(X^i)\right) \binom{n+d-2}{d-2}\\
& + (\text{polynomial of degree $<d-2$})
\end{align*}
for $n \gg 0$. It follows that 
\[
\rme_0(M)=j_0\rme_0(R/\fkp) \quad \text{and} \quad \rme_1(M)= j_0 \rme_0(R/\fkp)t + j_0 \rme_1(R/\fkp) - \sum_{i=0}^{j_0-1} \rme(X^i).
\]
Hence, $\rme_1(M)\le t\rme_0(M) + \ell_{R_\fkp}(M_\fkp){\cdot}\rme_1(R/\fkp)$ because $\rme(X^i)\ge 0$.

(b) (i) $\Rightarrow$ (ii): Suppose that $R/\fkp$ satisfies Serre's condition $(S_2)$. Because $\sum_{i=0}^{j_0-1} \rme(X^i)=0$, we obtain that $\dim_R X^i<d-1$ for all $0 \le i \le j_0-1$. Assume that $X^i\ne 0$, and choose $\fkq \in \Ass_R X^i$. Then, by localizing (\ref{*}) at $\fkq$, we obtain that $\depth_{R_\fkq} N^i_\fkq=0$.  Note that 
\[
\height_R \fkq = \height_{R/\fkp} (\fkq/\fkp) =\dim R/\fkp - \dim R/\fkq \ge 2
\] 
because $\dim_R X^i<d-1$ and $R/\fkp$ is catenary (see \cite[p.84]{matsu}). 
Hence, $\depth_{R_\fkq} (R/\fkp)_\fkq \ge 2$ because $R/\fkp$ satisfies $(S_2)$. By looking at \eqref{*}, it follows that $\depth_{R_\fkq} M^i_\fkq=0$, i.e., $\fkq \in \Ass_R M^i$. This contradicts  $\Ass_R M^i=\{\fkp \}$. Therefore, $X^i=0$ for all $0 \le i \le j_0-1$.

(b) (ii) $\Rightarrow$ (i): By (\ref{**}), we have 
\begin{align*} 
\ell_{R_0}(M_n)=&i_0{\cdot}\ell_{R_0}((R/\fkp)_{n-t})\\
=&i_0\rme_0(R/\fkp) \binom{n+d-1}{d-1} + \left( i_0 \rme_0(R/\fkp)(-t)  - i_0 \rme_1(R/\fkp) \right) \binom{n+d-2}{d-2}\\
&+ (\text{polynomial of degree $<d-2$})
\end{align*}
for $n\gg 0$. Hence, $\rme_0(M) = i_0\rme_0(R/\fkp)$ and $\rme_1(M) = i_0 t \rme_0(R/\fkp) + i_0 \rme_1(R/\fkp)$. Note that $\rme_0(M)=\ell_{R_\fkp}(M_\fkp){\cdot}\rme_0(R/\fkp)$ by (a); thus, $\rme_1(M)= t\rme_0(M) + \ell_{R_\fkp}(M_\fkp){\cdot}\rme_1(R/\fkp)$.
\end{proof}

\begin{corollary} 
Let $R=R_0[R_1]$ be a standard graded Noetherian domain of dimension $\ge 2$ such that $R_0$ is a field. Suppose that $R$ satisfies Serre's condition $(S_2)$. Let $M$ be a finitely generated graded torsionfree $R$-module such that $M=RM_t$. Then, 
$\rme_1(M)= t{\cdot}\rme_0(M) + \rme_1(R){\cdot}\rank_R M$ if and only if $M$ is an $R$-free module.
\end{corollary}

\begin{proof} 
We have only to show the ``only if" part. By noting that $R_0$ is a field, Theorem \ref{thm1} shows that there exists a family $\{M^i\}_{0\le i\le i_0}$ of graded residue modules of $M$ such that
\begin{align*}
\begin{split} 
0 \to R(-t) \to &M=M^0 \to M^1 \to 0,\\
0 \to R(-t) \to &M^1 \to M^2 \to 0,\\
&\vdots\\
0 \to R(-t) \to &M^{i_0-2} \to M^{i_0-1} \to 0, \text{ and}\\
0 \to R(-t) \to &M^{i_0-1} \to M^{i_0}=0 \to 0
\end{split}
\end{align*}
are graded exact sequences as $R$-modules. Since $R(-t)$ is a projective $R$-module, the exact sequence 
\[
0 \to R(-t) \to M^{i_0-2} \to M^{i_0-1}\cong R(-t) \to 0
\] 
splits, that is, $M^{i_0-2}  \cong R(-t)^2$. By descending induction on $i=0, \dots, i_0-2$, we obtain that $M\cong R(-t)^{i_0}$ and $i_0=\ell_{R_{(0)}}(M_{(0)})=\rank_R M$.
\end{proof}

\begin{ex} 
\begin{enumerate}[{\rm (a)}] 
\item Let $R=K[X, Y]$ be the polynomial ring over a field $K$ with $\deg X=\deg Y=1$. Set $M=R^{\oplus 2}/\left<\left(
\begin{smallmatrix}
X\\
Y 
\end{smallmatrix}
\right)\right>$. Then, $\Ass_R M=\{(0)\}$ and there exist graded exact sequences
\[
0 \to R \to M \to R/(X) \to 0 \quad \text{and} \quad 0 \to R \to M \to R/(Y) \to 0.
\]
In particular, graded commutative diagrams (\ref{*}) of $M$ are not unique.
\item Let $(R_0, \fkm_0)$ be an Artinian local ring, and let $R=R_0[X_1, \dots, X_d]$ be the polynomial ring over $R_0$ with $\deg X_j=1$ for $1 \le j \le d$. Set $i_0=\ell_{R_0}(R_0)$, and choose elements $0=f_0, f_1, \dots, f_{i_0}\in R_0$ such that $(f_0, \dots, f_{i+1})/(f_0, \dots, f_{i})\cong R_0/\fkm_0$ for $0 \le i \le i_0-1$. Then, a graded filtration (\ref{**}) of $R$ can be obtained as follows:
\begin{align*} 
0 \to (R/\fkm_0 R) \xrightarrow{f_1} &R \to R^1 \to 0,\\
0 \to (R/\fkm_0 R) \xrightarrow{f_2} &R^1 \to R^2 \to 0,\\
&\vdots\\
0 \to (R/\fkm_0 R) \xrightarrow{f_{i_0-1}} &R^{i_0-2} \to R^{i_0-1} \to 0, \text{ and}\\
0 \to (R/\fkm_0 R) \xrightarrow{f_{i_0}} &R^{i_0-1} \to R^{i_0}=0 \to 0.
\end{align*}
\end{enumerate} 
\end{ex}

\begin{proof} 
(a): Let $\overline{*}$ denote the image of $*\in R^2$ in $M$. It is straightforward to verify that the map $R\to M$, where $1\mapsto \overline{\left(
\begin{smallmatrix}
0\\
1 
\end{smallmatrix}
\right)}$, induces $0 \to R \to M \to R/(X) \to 0$. By symmetry, we also obtain that $0 \to R \to M \to R/(Y) \to 0$. It follows that $\Ass_R M=\{(0)\}$.

(b): Because $f_i$ is a socle of $R^{i-1}_0= A/(f_0, \dots, f_{i-1})$, we have $\fkm_0 R \subseteq (0):_R f_i$. This follows that $(0):_R f_i=\fkm_0 R$ because $\Ass_R R^{i-1} = \{\fkm_0 R\}$.
\end{proof}


\section{Application to Sally modules}\label{section3}

Throughout this section, let $(A, \fkm)$ be a Cohen--Macaulay local ring of dimension $d\ge 2$, and let  $I$ be an $\fkm$-primary ideal of $A$. Suppose that the residue field $A/\fkm$ is infinite. Then, there exists a parameter ideal $Q\subseteq I$ such that $Q$ is a reduction of $I$, that is, $I^{n+1}=QI^n$ for some $n \ge 0$ (\cite[Proposition 4.6.8]{BH}). 
Let 
\begin{align*} 
\calR(I)=A[It] \subseteq A[t] \quad \text{and} \quad  \calG(I)=\calR(I)/I\calR(I)\cong \bigoplus_{n\ge 0}I^n/I^{n+1},
\end{align*}
be the {\it Rees algebra} of $I$ and the {\it associated graded ring} of $I$, respectively, where $A[t]$ is the polynomial ring over $A$. To investigate the Hilbert coefficients of $I$ relating to the structures of $\calR(I)$ and $\calG(I)$, we further define the finitely generated graded $\calR(Q)$-module
\begin{align*} 
&S=\calS_Q (I)=I \calR(I)/I \calR(Q)\cong \bigoplus_{n\ge 1} I^{n+1}/Q^{n}I.
\end{align*}
$S$ is called the {\it Sally module} of $I$ with respect to $Q$ (\cite{V}). 
The importance of the notion of Sally modules is to describe a correction term of the Hilbert function of $I$. We summarize the properties of Sally modules below.

\begin{proposition}{\rm (\cite[Lemma 2.1 and Proposition 2.2]{GNO})}\label{basic}
The following hold true:
\begin{enumerate}[{\rm (a)}] 
\item We have  
\begin{align*} 
\ell_A(A/I^{n+1})=&\rme_0(I)\binom{n+d}{d}-[\rme_0(I)-\ell_A (A/I)]\binom{n+d-1}{d-1}- \ell_A(S_n) 
\end{align*}
for all $n\ge 0$.
\item Let $\calM=\fkm \calR(Q)+\calR(Q)_+$ denote the graded maximal ideal of $\calR(Q)$. Then, for $n>0$, $(S/\calM S)_n=0$ if and only if $I^{n+1}=QI^{n}$. 
\item $\fkm^\ell S=0$ for $\ell \gg 0$.
\item $\Ass_{\calR(Q)}S\subseteq \{\fkm \calR(Q)\}$. In particular, either $S=0$ or $\dim_{\calR(Q)} S=d$ holds.
\item $\depth \calG(I)\ge d-1$ if and only if either $S=0$ or $S$ is a ($d$-dimensional) Cohen--Macaulay $\calR(Q)$-module. If $S$ is neither $0$ nor a Cohen--Macaulay $\calR(Q)$-module, then $\depth \calG(I)= \depth_{\calR(Q)} S-1$.
\end{enumerate} 
\end{proposition}

The next proposition, Proposition \ref{prop3.2}, is based on the following observation: we can regard $S$ as a graded $V=\calR(Q)/\fkm^\ell \calR(Q)$-module for $\ell \gg 0$ by Proposition \ref{basic}(c). Here, $V$ is a standard graded Noetherian ring such that $V_0=A/\fkm^\ell$ is an Artinian local ring. Since $V/\fkm V \cong \calR(Q)/\fkm \calR(Q) \cong (A/\fkm)[X_1, \dots, X_d]$ is a polynomial ring over the field $A/\fkm$ (\cite[Theorem 1.1.8]{BH}), $\fkm V$ is a prime ideal of $V$. Hence, according to Proposition \ref{basic}(d), the associated prime ideals of $S$ (as a $V$-module) are only $\fkm V$. Therefore, we can apply Theorem \ref{thm1} to submodules of $S$ generated in single degree. Let 
\[
r=\mathrm{red}_Q I=\min\{n\ge 0 \mid I^{n+1}=QI^n\}
\] 
denote the {\it reduction number} of $I$ with respect to $Q$. The graded minimal generators of $S$ are in degree $1$ to $r-1$ by Proposition \ref{basic}(b). 

\begin{proposition} \label{prop3.2}
Let $(A, \fkm)$ be a Cohen--Macaulay local ring of dimension $d\ge 2$, and let  $I$ be an $\fkm$-primary ideal of $A$. Suppose that $A/\fkm$ is an infinite field, and choose a parameter ideal $Q\subseteq I$ of $A$ as a reduction of $I$. 
If $r = \mathrm{red}_Q I \ge 2$, then the following are true:
\begin{enumerate}[{\rm (a)}] 
\item $\ell_A(A/I)\ge \rme_0(I) - \rme_1(I) + \dfrac{\rme_2(I)}{r-1}$ holds.  
\item $\ell_A(A/I)= \rme_0(I) - \rme_1(I) + \dfrac{\rme_2(I)}{r-1}$ implies that $r=2$.
\end{enumerate} 
\end{proposition}

\begin{proof} 
(a): Let $S'=\calR(Q) S_{r-1}$ denote the submodule of $S$ generated in degree $r-1$. Then, because $S'_{n}=S_n$ for all $n\ge r-1$, we obtain that 
\begin{align*} 
\ell_A(A/I^{n+1}) = & \rme_0(I)\binom{n+d}{d}-[\rme_0(I)-\ell_A (A/I)]\binom{n+d-1}{d-1}- \ell_A(S'_n) \\
=& \rme_0(I)\binom{n+d}{d}-[\rme_0(I)-\ell_A (A/I) + \rme_0(S')]\binom{n+d-1}{d-1} \\
 & + \rme_1(S')\binom{n+d-2}{d-2} +\cdots +(-1)^d \rme_{d-1}(S')
\end{align*}
for $n \gg 0$. Hence, $\rme_1(I)=\rme_0(I)-\ell_A (A/I) + \rme_0(S')$ and $\rme_2(I)=\rme_1(S')$. By Theorem \ref{thm1}(a) (see the observation before Proposition \ref{prop3.2}), we have 
\[
\rme_1(S') \le (r-1)\rme_0(S') + \ell_{\calR(Q)_\fkp}(S'_\fkp){\cdot}\rme_1(\calR(Q)/\fkp) = (r-1)\rme_0(S') 
\]  
since $\calR(Q)/\fkp = \calR(Q)/\fkm \calR(Q) \cong (A/\fkm)[X_1, X_2, \dots, X_d]$ and thus $\rme_1(\calR(Q)/\fkp)=0$. Hence, $\rme_1(I) \ge \rme_0(I)-\ell_A (A/I) +\dfrac{\rme_2(I)}{r-1}$.

(b): Note that $\calR(Q)/\fkp$ is a Cohen--Macaulay ring of dimension $\ge 2$ and thus satisfies Serre's condition $(S_2)$. Hence, the equality $\ell_A(A/I)= \rme_0(I) - \rme_1(I) + \dfrac{\rme_2(I)}{r-1}$ holds, which is equivalent to saying that $\rme_1(S') = (r-1)\rme_0(S')$ holds, if and only if $S'$ has a graded filtration (\ref{**}) by Theorem \ref{thm1}(b). It follows that $S'$ is a $d$-dimensional Cohen--Macaulay $\calR(Q)$-module. Indeed, we apply the depth formula to each exact sequence in the graded filtration (\ref{**})  of $S'$recursively from $i=i_0-1$ to $0$, in that order. Then, we obtain that all residue modules $S'^{\,i}$ of $S'$ appearing in \eqref{**} are of depth $\ge d$. On the other hand, $S'^{\, i}$ is of dimension $d$ since $\Ass_{\calR(Q)} S'^{\, i}=\{\fkp\}$. Hence, $S'^{\,i}$ are $d$-dimensional Cohen--Macaulay $\calR(Q)$-modules. In particular, $S'=S'^{\, 0}$ is a $d$-dimensional Cohen--Macaulay $\calR(Q)$-module. 

We prove that $S=S'$. Indeed, if $S\ne S'$, by applying the depth formula to the exact sequence $0 \to S' \to S \to S/S' \to 0$, we obtain that $\depth_{\calR(Q)} S=0$  since $S/S'$ is of finite length. This contradicts $\Ass_{\calR(Q)}S= \{\fkp\}$. Therefore, $S=S'$. By noting that $S_0=0$ and $S_1\cong I^2/QI\ne 0$ because of the assumption that $r\ge 2$, the fact that $S=S'$ is generated in single degree $r-1$ forces $r=2$.
\end{proof}

By focusing on the case of $r=2$, we obtain the second main result of this study.

\begin{theorem} \label{thm3.3}
Suppose that $I^3=QI^2$. Then, $\ell_A(A/I)\ge \rme_0(I) - \rme_1(I) + \rme_2(I)$ holds. Furthermore, the following are equivalent:
\begin{enumerate}[{\rm (a)}] 
\item $\ell_A(A/I)= \rme_0(I) - \rme_1(I) + \rme_2(I)$,
\item $S$ has a filtration {\rm (\ref{**})} (we include the case that $S=0$), and 
\item $\depth \calG(I) \ge d-1$.
\end{enumerate} 
When the equivalence conditions hold, 
\begin{center} 
$\ell_A (A/I^{n+1})=\rme_0(I) \binom{n+d}{d} -\rme_1(I)\binom{n+d-1}{d-1}+\rme_2(I)\binom{n+d-2}{d-2}$ \quad for all $n\ge 0$.
\end{center}
\end{theorem}

\begin{proof}
Suppose that $I^2=QI$. Then, $\ell_A(A/I)= \rme_0(I) - \rme_1(I)$ and $\rme_2(I)=0$ by \cite[Theorem 2.1]{H}. $S=0$ by Proposition \ref{basic}(b). $\calG(I)$ is a Cohen--Macaulay ring of dimension $d$ by \cite[Proposition 3.1]{VV}. Thus, (a), (b), and (c) hold. Then, the Hilbert function is also calculated as $\ell_A (A/I^{n+1})=\rme_0(I) \binom{n+d}{d} -\rme_1(I)\binom{n+d-1}{d-1}$ for all $n\ge 0$ by \cite[Theorem 2.1]{H}. Hence, we may assume that $I^2\ne QI$, that is, $r=2$. 
Then, the inequality $\ell_A(A/I)\ge \rme_0(I) - \rme_1(I) + \rme_2(I)$ follows from Proposition \ref{prop3.2}. We prove the equivalence of (a), (b), and (c).

(a) $\Rightarrow$ (b): By recalling the proof of Proposition \ref{prop3.2}(a), $\ell_A(A/I)= \rme_0(I) - \rme_1(I) + \rme_2(I)$ holds if and only if $S'=\calR(Q) S_{1}$ has a graded filtration (\ref{**}). In addition, since $r=2$, we have $S=S'$ by Proposition \ref{basic}(b).

(b) $\Rightarrow$ (c): If $S$ has a filtration {\rm (\ref{**})}, then $S$ is a $d$-dimensional Cohen--Macaulay $\calR(Q)$-module by applying the depth formula to each exact sequence in the filtration {\rm (\ref{**})} of $S$ recursively from the very bottom. It follows that $\depth \calG(I) \ge d-1$ by Proposition \ref{basic}(e). 

(c) $\Rightarrow$ (a): If $\depth \calG(I) \ge d-1$, then we have that 
\begin{align*} 
\rme_1(I)=\sum_{i=0}^{r-1}\ell_A(I^{i+1}/QI^i) \quad \text{and} \quad 
\rme_2(I)=\sum_{i=1}^{r-1} i{\cdot}\ell_A(I^{i+1}/QI^i)
\end{align*}
by \cite[Theorem 4.7(b)]{HM} and \cite[Theorem 3.1]{CPR}, respectively (or see \cite[Theorem 2.5]{RV}). 
Because we suppose that $r=2$, it follows that 
\begin{align*} 
\rme_1(I)=\ell_A(I/Q) + \ell_A(I^2/QI) \quad \text{and} \quad 
\rme_2(I)=\ell_A(I^2/QI).
\end{align*}
Thus, $\rme_1(I)=\ell_A(A/Q) - \ell_A(A/I) + \rme_2(I)$; hence, $\ell_A(A/I)= \rme_0(I) - \rme_1(I) + \rme_2(I)$.

When the equivalence conditions hold, since $S$ has a filtration {\rm (\ref{**})}, Proposition \ref{basic}(a) proves that
\begin{align*} 
\ell_A (A/I^{n+1})=& \rme_0(I)\binom{n+d}{d}-[\rme_0(I)-\ell_A (A/I)]\binom{n+d-1}{d-1} - \ell_A(S_1){\cdot}\ell_A([\calR(Q)/\fkp]_{n-1}) \\
=& \rme_0(I)\binom{n+d}{d}-[\rme_0(I)-\ell_A (A/I)]\binom{n+d-1}{d-1} - \ell_A(I^2/QI)\binom{n-1+d-1}{d-1} \\
=&\rme_0(I) \binom{n+d}{d} -\rme_1(I)\binom{n+d-1}{d-1}+\rme_2(I)\binom{n+d-2}{d-2}
\end{align*}
for all $n\ge 0$.
\end{proof}

Note that Theorem \ref{thm3.3} improves \cite[Theorem 3.3]{K}.
Furthermore, in contrast to Theorem \ref{thm3.3}, it is known that the converse inequality $\ell_A(A/I)\le \rme_0(I) - \rme_1(I) + \rme_2(I)$ holds if $I$ is either an integrally closed ideal (\cite[Theorem 12]{I}) or a Ratliff--Rush closed ideal of height two (\cite[(2.5) Corollary]{S2}). Hence, by combining Theorem \ref{thm3.3} and \cite{I, S2}, we obtain the following.

\begin{corollary} {\rm (cf. \cite[Theorem 3.6]{CPR})} \label{cor3,4}
Suppose that $I^3=QI^2$ and either 
\begin{enumerate}[{\rm (a)}] 
\item $I$ is an integrally closed ideal or
\item $d=2$ and $I$ is a Ratliff--Rush closed ideal in the sense of \cite[Chapter ${\rm V\hspace{-.15em}I\hspace{-.15em}I\hspace{-.15em}I}$, Notation]{Mc}.
\end{enumerate}
Then, $\depth \calG(I) \ge d-1$ and 
$\ell_A (A/I^{n+1})=\rme_0(I) \binom{n+d}{d} -\rme_1(I)\binom{n+d-1}{d-1}+\rme_2(I)\binom{n+d-2}{d-2}$ \quad for all $n\ge 0$.
\end{corollary}

The following lemmas are useful for providing examples of Theorem \ref{thm3.3}. Let $\ol{*}$ denote the {\it integral closure} of an ideal $*$.

\begin{lemma} \label{lem3.5}
Assume that $I^2=\ol{I}^2=Q\ol{I}$. Then, $I^3=QI^2$ and $\ell_A(I^2/QI)=d{\cdot}\ell_A(\ol{I}/I)$. In particular, $\mathrm{red}_Q I=2$ if $I\ne \ol{I}$.
\end{lemma}

\begin{proof} 
By assumption, we obtain that $I^3=I{\cdot}\ol{I}^2=I{\cdot}Q\ol{I}=QI^2$. The latter equality follows from 
\[
\ell_A(I^2/QI)=\ell_A(Q\ol{I}/QI)=\ell_A(Q/QI)- \ell_A(Q/Q\ol{I})=d{\cdot}(\ell_A(A/I)-\ell_A(A/\ol{I})),
\]
where the last equality follows from the isomorphisms $Q/QI\cong (Q/Q^2)\otimes_A A/I\cong (A/I)^{\oplus d}$ and $Q/Q\ol{I} \cong (A/\ol{I})^{\oplus d}$.
\end{proof}

\begin{lemma} \label{lem3.6}
\begin{enumerate}[{\rm (a)}] 
\item Let $R=K[X, Y]$ be the polynomial ring over a field $K$ with $\deg X=\deg Y =1$. Let $I_0$ be a monomial ideal generated in single degree $t$, and let $Q_0=(X^t, Y^t)$. Suppose that $Q_0\subseteq I_0$. Then, $\ol{Q_0}=\ol{I_0}=(X, Y)^t$ and $\mathrm{red}_{Q_0} \ol{I_0}= 1$.
\item Let $A=K[[X, Y]]$ be the formal power series ring over a field $K$ with the maximal ideal $\fkm=(X, Y)$. Let $0=n_1<n_2<\cdots<n_s= t$ be integers, and let $I$ be an ideal generated by
\[
Y^t=X^{n_1}Y^{t-n_1}, X^{n_2}Y^{t-n_2}, \dots, X^{n_s}Y^{t-n_s}=X^t.
\]
Then, $Q=(X^t, Y^t)$ is a reduction of $I$. Furthermore, if $Q \subsetneq I \subsetneq \fkm^t$ and $I^2=\fkm^{2t}$, then $\mathrm{red}_Q I=2$.
\end{enumerate} 
\end{lemma}

\begin{proof} 
(a): To prove the former equality, because $Q_0\subseteq I_0 \subseteq (X, Y)^t$, we need only to demonstrate that $(X, Y)^t\subseteq \ol{Q_0}$.
Indeed, this follows from the fact that $(X^s Y^{t-s})^t\in Q_0^t$ for all $0 \le s \le t$. The latter inequality follows from $\ol{I_0}^2=(X, Y)^{2t}=Q_0(X, Y)^t=Q_0\ol{I_0}$.

(b): This follows from (a) and Lemma \ref{lem3.5}.
\end{proof}

\begin{ex} 
Let $A=K[[X, Y]]$ be the formal power series ring over a field $K$. Set $Q=(X^7, Y^7)$ and $I=Q + (X^6Y, X^5Y^2, X^2Y^5, XY^6)$. Then, $\mathrm{red}_Q I=2$ and 
\begin{align*} 
\ell_A(A/I^{n+1})=\begin{cases} 
31 & (n=0)\\
49\binom{n+2}{2}-21\binom{n+1}{1} & (n\ge 1).
\end{cases}
\end{align*}
It follows that $\ell_A(A/I)=31 > \rme_0(I)-\rme_1(I)+\rme_2(I)=28$; hence, $\depth \calG(I)=0$.
\end{ex}

\begin{proof} 
Set $\fkm=(X, Y)$.
It is routine to prove that $I^2=\fkm^{14}$, and this demonstrates that $\mathrm{red}_Q I=2$ by Lemma \ref{lem3.6}.
It is thus proven that $\ell_A(A/I^{n+1})=\ell_A(A/\fkm^{7n+7})=49\binom{n+2}{2}-21\binom{n+1}{1}$ for all $n \ge 1$. 
\end{proof}

\begin{ex} 
Let $s>0$, and let $A=K[[XY^i \mid 0 \le i \le 2s+1]]$ be a subring of the formal power series ring $K[[X, Y]]$ over a field $K$. Then, $A$ is a Cohen--Macaulay local normal domain of dimension $2$  (for example, \cite[Theorem 6.3.5]{BH}). Set $I=(XY^i \mid 0 \le i \le s) + (XY^{2s+1})$ and $Q=(X, XY^{2s+1})$. Then, the following hold true:
\begin{enumerate}[{\rm (a)}] 
\item $I^2\subseteq Q$ and $\ell_A(I^2/QI)=s$. In particular, $\calG(I)$ is not a Cohen--Macaulay ring.
\item $I^3=QI^2$.
\item $\ell_A(A/I^{n+1})=(2s+1)\binom{n+2}{2} - 2s\binom{n+1}{1} + s$ for all $n\ge 0$.
\item $\depth \calG(I)=1$.
\end{enumerate} 
\end{ex}

\begin{proof} 
Let $\fkm$ denote the maximal ideal of $A$.

(a): It is easy to verify that $I^2 \subseteq \fkm^2\subseteq Q$. It follows that $\fkm I^2\subseteq \fkm^2 I \subseteq QI \subseteq I^2$. Hence, 
\[
I^2=(X^2Y^i \mid 0 \le i \le 3s+1 \text{ or } i=4s+2) = QI+(X^2Y^i \mid s+1 \le i \le 2s)
\] 
demonstrates that $\ell_A(I^2/QI)=s$. Thus, $I^2\cap Q=I^2\ne QI$; hence, $\calG(I)$ is not a Cohen--Macaulay ring by \cite{VV}.

(b): This is straightforward. 

(c): By induction on $n$, we obtain 
\[
I^n=(X^nY^i \mid \text{$0 \le i \le (2s+1)(n-1)+s$ \quad or \quad $i=(2s+1)n$}).
\]
Hence, 
\begin{align*} 
\ell_A(A/I^{n+1}) =& \ell_A(A/\fkm^{n+1}) + s = \sum_{k=0}^{n} [(2s+1)k + 1] + s\\
=& (2s+1)\binom{n+2}{2} - 2s\binom{n+1}{1} + s
\end{align*}
for all $n\ge 0$.

(d): Because $\ell_A(A/I^{n+1})$ agrees with the polynomial of (c) at $n=0$, $\depth \calG(I)>0$ according to Theorem \ref{thm3.3}. Therefore, $\depth \calG(I)=1$ by (a).
\end{proof}

In the remainder of this paper, we present the case in which $I$ is integrally closed. Set 
\begin{align*} 
&L=\calL_Q (I)= \calR(Q)S_1\cong \bigoplus_{n\ge 1} Q^{n-1}I^{2}/Q^nI\quad \text{and}\\
&C=\calC_Q (I)=(I^2 \calR(I)/I^2 \calR(Q))(-1)\cong \bigoplus_{n\ge 2} I^{n+1}/Q^{n-1}I^2.
\end{align*}
Then, we have the graded exact sequence
\[
0 \to L \to S \to C \to 0
\]
of $\calR(Q)$-modules. This is certainly part of the filtration of the Sally module defined by Vaz Pinto \cite{VP}.  This filtration is effective if $I^2=Q\cap I$ as follows:

\begin{proposition}\label{basic2} {\rm (\cite[Propositions 2.2 and 2.8, and Lemma 2.11]{OR})}
Suppose that $Q\cap I^2=QI$ holds. Then, we have the following.
\begin{enumerate}[{\rm (a)}] 
\item We have
\begin{align*} 
\ell_A(A/I^{n+1})=&\rme_0(I)\binom{n+d}{d}-[\rme_0(I)-\ell_A (A/I)+\ell_A(I^2/QI)]\binom{n+d-1}{d-1} \\
&+\ell_A(I^2/QI)\binom{n+d-2}{d-2}- \ell_A(C_n) 
\end{align*}
for all $n\ge 0$.
\item $\Ass_{\calR(Q)}C\subseteq \{\fkm \calR(Q)\}$. In particular, either $C=0$ or $\dim_{\calR(Q)} C=d$ holds.
\item $\depth \calG(I)\ge d-1$ if and only if either $C=0$ holds or $C$ is a ($d$-dimensional) Cohen--Macaulay $\calR(Q)$-module. If $C$ is neither $0$ nor a Cohen--Macaulay $\calR(Q)$-module, then $\depth \calG(I)= \depth_{\calR(Q)} C-1$.
\end{enumerate} 
\end{proposition}

\begin{rem} \label{rem1} {\rm (\cite[Theorem 1]{I2})}
If $I$ is an integrally closed ideal, then $Q\cap I^2=QI$ holds.
\end{rem}

By using Proposition \ref{basic2}, we can apply Theorem \ref{thm1} to $C$ by regarding $C$ as a graded $V=\calR(Q)/\fkm^\ell \calR(Q)$-module (recall the observation before Proposition \ref{prop3.2}). The proofs of Proposition \ref{finalprop} and Theorem \ref{final} proceed in the same manner as in Proposition \ref{prop3.2} and Theorem \ref{thm3.3}. We include the proofs for the sake of completeness.

\begin{proposition} \label{finalprop}
Suppose that $I$ is integrally closed and $r =\mathrm{red}_Q I\ge 3$. Then, the following are true:
\begin{enumerate}[{\rm (a)}] 
\item $\ell_A(A/I)\ge \rme_0(I) - \rme_1(I) + \frac{(r-2)\ell_A(I^2/QI) + \rme_2(I)}{r-1}$ holds.  
\item $\ell_A(A/I)= \rme_0(I) - \rme_1(I) + \frac{(r-2)\ell_A(I^2/QI) + \rme_2(I)}{r-1}$ implies that $r=3$.
\end{enumerate} 
\end{proposition}


\begin{proof} 
(a): Let $C'=\calR(Q) C_{r-1}$ denote the submodule of $C$ generated in degree $r-1$. Then, because $C'_{n}=C_n$ for all $n\ge r-1$, we obtain that 
\begin{align*} 
\ell_A(A/I^{n+1}) = & \rme_0(I)\binom{n+d}{d}-[\rme_0(I)-\ell_A (A/I)+\ell_A(I^2/QI)]\binom{n+d-1}{d-1}\\
&+\ell_A(I^2/QI)\binom{n+d-2}{d-2} - \ell_A(C'_n) \\
=& \rme_0(I)\binom{n+d}{d}-[\rme_0(I)-\ell_A (A/I) +\ell_A(I^2/QI) + \rme_0(C')]\binom{n+d-1}{d-1} \\
 & +[\ell_A(I^2/QI) + \rme_1(C')]\binom{n+d-2}{d-2} - \rme_2(C')\binom{n+d-3}{d-3} +\cdots +(-1)^d \rme_{d-1}(C')
\end{align*}
for $n \gg 0$. Hence, $\rme_1(I)=\rme_0(I)-\ell_A (A/I) +\ell_A(I^2/QI) + \rme_0(C')$ and $\rme_2(I)=\ell_A(I^2/QI) + \rme_1(C')$. By Theorem \ref{thm1}(a), we have 
\[
\rme_1(C') \le (r-1)\rme_0(C') + \ell_{\calR(Q)_\fkp}(C'_\fkp){\cdot}\rme_1(\calR(Q)/\fkp) = (r-1)\rme_0(C'). 
\]  
Hence, $\rme_1(I) \ge \rme_0(I)-\ell_A (A/I) +\ell_A(I^2/QI) +\dfrac{\rme_2(I) - \ell_A(I^2/QI)}{r-1}$.

(b): Note that $\calR(Q)/\fkp$ is a Cohen--Macaulay ring of dimension $\ge 2$ and thus satisfies Serre's condition $(S_2)$. Hence, the equality $\ell_A(A/I)= \rme_0(I) - \rme_1(I) + \frac{(r-2)\ell_A(I^2/QI) + \rme_2(I)}{r-1}$ holds if and only if $C'$ has a graded filtration (\ref{**}) by Theorem \ref{thm1}(b). It follows that $C'$ is a $d$-dimensional Cohen--Macaulay $\calR(Q)$-module. Indeed, we apply the depth formula to each exact sequence in the graded filtration (\ref{**})  of $C'$ recursively from $i=i_0-1$ to $0$, in that order. Then, we obtain that all residue modules $C'^{\,i}$ of $C'$ appearing in \eqref{**} are of depth $\ge d$. On the other hand, $C'^{\, i}$ is of dimension $d$ since $\Ass_{\calR(Q)} C'^{\, i}=\{\fkp\}$. Hence, $C'^{\,i}$ are $d$-dimensional Cohen--Macaulay $\calR(Q)$-modules. In particular, $C'=C'^{\, 0}$ is a $d$-dimensional Cohen--Macaulay $\calR(Q)$-module. 

We prove that $C=C'$. Indeed, if $C\ne C'$, by applying the depth formula to the exact sequence $0 \to C' \to C \to C/C' \to 0$, we obtain that $\depth_{\calR(Q)} C=0$ since $C/C'$ is of finite length. This contradicts $\Ass_{\calR(Q)}C= \{\fkp\}$. Therefore, $C=C'$. By noting that $C_1=0$ and $C_2\cong I^3/QI^2\ne 0$ because of the assumption that $r\ge 3$, the fact that $C=C'$ is generated in single degree $r-1$ forces $r=3$.
\end{proof}


\begin{theorem} \label{final}
Suppose that $I$ is integrally closed and  $I^4=QI^3$. Then, $\ell_A(A/I)\ge \rme_0(I) - \rme_1(I) + \frac{\ell_A(I^2/QI) + \rme_2(I)}{2}$ holds. Furthermore, the following are equivalent:
\begin{enumerate}[{\rm (a)}] 
\item $\ell_A(A/I)= \rme_0(I) - \rme_1(I) + \frac{\ell_A(I^2/QI) + \rme_2(I)}{2}$,
\item $C$ has a filtration {\rm (\ref{**})} (we include the case that $C=0$), and 
\item $\depth \calG(I) \ge d-1$.
\end{enumerate} 
When this is the case, 
\[
\ell_A (A/I^{n+1})=\rme_0(I) \binom{n+d}{d} -\rme_1(I)\binom{n+d-1}{d-1}+\rme_2(I)\binom{n+d-2}{d-2} - \rme_3(I)\binom{n+d-3}{d-3}
\] 
for all $n\ge 1$ (If $d=2$, we consider that $\binom{n+d-3}{d-3}=0$).
\end{theorem}

\begin{proof}
Suppose that $I^3=QI^2$. Then, $S$ is generated in degree $1$, thus $C=S/L=0$. We obtain that $\ell_A(A/I)= \rme_0(I) - \rme_1(I)+\rme_2(I)$ and $\depth \calG(I) \ge d-1$ by Corollary \ref{cor3,4}. By recalling that $\rme_2(I)=\ell_A(I^2/QI)$ (see the proof of Theorem \ref{thm3.3}(c)$\Rightarrow$(a)), the assertions (a), (b), and (c) hold. The Hilbert function is also calculated as desired by Corollary \ref{cor3,4}. Hence, we may assume that $I^3\ne QI^2$, that is, $r=3$. 
Then, the inequality $\ell_A(A/I)\ge \rme_0(I) - \rme_1(I) + \frac{\ell_A(I^2/QI) + \rme_2(I)}{2}$ follows from Proposition \ref{finalprop}(a). We prove the equivalence of (a), (b), and (c).

(a) $\Rightarrow$ (b): By recalling the proof of Proposition \ref{finalprop}(a), $\ell_A(A/I)= \rme_0(I) - \rme_1(I) + \frac{\ell_A(I^2/QI) + \rme_2(I)}{2}$ holds if and only if $C'=\calR(Q) C_{2}$ has a graded filtration (\ref{**}). Furthermore, the proof of Proposition \ref{finalprop}(b) also shows that $C=C'$ when $C'$ has a graded filtration (\ref{**}). 

(b) $\Rightarrow$ (c): If $C$ has a filtration {\rm (\ref{**})}, then $C$ is a $d$-dimensional Cohen--Macaulay $\calR(Q)$-module by applying the depth formula to each exact sequence in the filtration {\rm (\ref{**})} of $C$ recursively from the very bottom. It follows that $\depth \calG(I) \ge d-1$ by Proposition \ref{basic2}(c).

(c) $\Rightarrow$ (a): If $\depth \calG(I) \ge d-1$, then we have that 
\begin{align*} 
\rme_1(I)=\sum_{i=0}^{r-1}\ell_A(I^{i+1}/QI^i) \quad \text{and} \quad 
\rme_2(I)=\sum_{i=1}^{r-1} i{\cdot}\ell_A(I^{i+1}/QI^i)
\end{align*}
by \cite[Theorem 4.7(b)]{HM} and \cite[Theorem 3.1]{CPR}, respectively. 
Because we suppose that $r=3$, it follows that 
\begin{align*} 
\rme_1(I)=\ell_A(I/Q) + \ell_A(I^2/QI) + \ell_A(I^3/QI^2) \quad \text{and} \quad 
\rme_2(I)=\ell_A(I^2/QI) + 2\ell_A(I^3/QI^2).
\end{align*}
Thus, $\rme_1(I)=\ell_A(A/Q) - \ell_A(A/I) +  \ell_A(I^2/QI) + \frac{\rme_2(I) -  \ell_A(I^2/QI)}{2}$; hence, $\ell_A(A/I)= \rme_0(I) - \rme_1(I) + \frac{\ell_A(I^2/QI) + \rme_2(I)}{2}$.

When the equivalence conditions hold, since $C$ has a filtration {\rm (\ref{**})}, Proposition \ref{basic2}(a) proves that
\begin{align*} 
\ell_A(A/I^{n+1})=&\rme_0(I)\binom{n+d}{d}-[\rme_0(I)-\ell_A (A/I)+\ell_A(I^2/QI)]\binom{n+d-1}{d-1}\\
&+\ell_A(I^2/QI)\binom{n+d-2}{d-2}- \ell_A(C_2){\cdot}\ell_A([\calR(Q)/\fkp]_{n-2})\\
=&\rme_0(I)\binom{n+d}{d}-[\rme_0(I)-\ell_A (A/I)+\ell_A(I^2/QI)]\binom{n+d-1}{d-1}\\
&+\ell_A(I^2/QI)\binom{n+d-2}{d-2}- \ell_A(I^3/QI^2)\binom{n-2+d-1}{d-1} 
\end{align*}
for all $n\ge 0$. By noting that 
\[
\binom{n-2+d-1}{d-1} = \binom{n+d-1}{d-1} -2 \binom{n+d-2}{d-2} + \binom{n+d-3}{d-3}
\]
for all $n\ge 0$, we obtain the conclusion as desired.
\end{proof}

\begin{ex}{\rm (\cite[Theorem 5.2]{OR})}
Let $m>0$, $d\ge 2$, and $D=K[[\{X_j\}_{1\le j \le m},  \{Y_i\}_{1\le i\le d}, \{Z_i\}_{1\le i\le d}]]$ be the formal power series  ring over a field $K$. Let
\[
\fka=(X_1, \dots, X_{m}){\cdot}(X_1, \dots, X_{m}, Y_1, \dots, Y_d) + (Y_iY_j \mid 1\le i, j\le d, i\ne j) + (Y_i^3-Z_iX_{m} \mid 1\le i \le d).
\]
Set $A=D/\fka$.  We denote by $x_*, y_*, z_*$ the image of $X_*, Y_*, Z_*$ into $A$. Let $\fkm$ be the maximal ideal of $A$ and 
$Q=(z_1, \dots, z_d)$. Then, the following assertions are true:
\begin{enumerate}[{\rm (a)}] 
\item $A$ is a Cohen--Macaulay local ring of dimension $d$.
\item $\fkm^4=Q\fkm^3$.
\item $\rme_0(\fkm)=m+2d+1, \rme_1(\fkm)=m+3d+1, \rme_2(\fkm)=d+1, \rme_i(\fkm)=0$ for all $3\le i \le d$, and  $\ell_A(\fkm^2/Q\fkm)=d$. Hence, 
\[
\ell_A(A/\fkm)=1> \rme_0(\fkm) - \rme_1(\fkm) + \frac{\ell_A(\fkm^2/Q\fkm) + \rme_2(\fkm)}{2}=\frac{1}{2}.
\]
\item $\depth \calG(\fkm) = 0$.
\end{enumerate}  
\end{ex}

\begin{acknowledgments}
The author was supported by JSPS KAKENHI Grant Numbers JP19J10579 and JP21K13766. The author is grateful to Kazuho Ozeki and Ryotaro Isobe for giving helpful comments. In particular, Kazuho Ozeki improved Theorem \ref{thm1}(b)(ii) $\Rightarrow$ (i) and Ryotaro Isobe gave the proof of Theorem \ref{thm3.3}(c) $\Rightarrow (a)$. The author also would like to thank the referee for his/her valuable comments. The proof of Lemma \ref{2.4} is improved by the referee.
\end{acknowledgments}


\end{document}